\documentclass[oneside,a4paper,reqno,11pt]{amsart}
\usepackage{latexsym,amsfonts,amssymb,amsmath,amscd}
\usepackage{graphicx,color,ifpdf}
\allowdisplaybreaks

\usepackage[bookmarks=true,pdfstartview=FitH, pdfborder={0 0 0},
colorlinks=true,citecolor=blue, linkcolor=blue]{hyperref}

\usepackage{aliascnt}

\newtheorem*{ques}{Question 1}
\newtheorem*{prop}{Proposition 2}
\newtheorem*{cor5}{Corollary 5}
\newtheorem*{cor6}{Corollary 6}

\theoremstyle{definition}
\newtheorem*{rem3}{Remark 3}
\newtheorem*{rem4}{Remark 4}

\theoremstyle{theorem}

\newaliascnt{lemma}{theorem}

\aliascntresetthe{lemma}

\newaliascnt{proposition}{theorem}

\aliascntresetthe{proposition}

\newaliascnt{corollary}{theorem}

\aliascntresetthe{corollary}

\theoremstyle{definition}
\newaliascnt{definition}{theorem}

\aliascntresetthe{definition}

\newaliascnt{example}{theorem}

\aliascntresetthe{example}

\newaliascnt{exercise}{theorem}

\aliascntresetthe{exercise}

\newaliascnt{question}{theorem}

\aliascntresetthe{question}

\newaliascnt{problem}{theorem}

\aliascntresetthe{problem}

\theoremstyle{remark}
\newaliascnt{remark}{theorem}

\aliascntresetthe{remark}

\newaliascnt{notation}{theorem}

\aliascntresetthe{notation}

\newaliascnt{fact}{theorem}

\aliascntresetthe{fact}

\numberwithin{equation}{theorem}%
\numberwithin{figure}{theorem}

\renewcommand{\Bbb}[1]{\mathbb{#1}}

\begin{document}
\title{an injectivity radius estimate in terms of metric 
sphere}
\author{Shicheng Xu}
\email{shichengxu@gmail.com}
\address{School of Mathematical Sciences, Capital Normal University, 
Beijing, China}

\thanks{\it 2000 Mathematics Subject Classification.\rm\ 53C20. 53C35}
\thanks{Project 11171143 supported by NSFC}
\thanks{Keywords: Injectivity radius, Complete manifold}

\date{\today}

\begin{abstract}
In this paper we prove that if a point $p$ in a complete 
Riemannian manifold is not a cut point of any 
point whose distance to $p$ is $r$, then the 
injectivity radius of $p$ is strictly large than $r$. As 
a corollary we give a positive answer to a problem 
raised by Z. Sun and J. Wan.
\end{abstract}
\maketitle

This paper is to answer a question asked by Z. Sun and J. 
Wan in \cite{SW14}. Let $M$ be a complete noncompact 
Riemannian manifold, and let $i_p$ denote the injectivity 
radius at $p$ of $M$. Let $$i(p,r)=\min \{i_x : 
\forall\, x\in M \text{ s.t. }\operatorname{d}(x,p)= r\},$$ 
where $\operatorname{d}(x,p)$ is the distance between two 
points $x$ and $p$.
According to \cite{SW14}, they defined a number 
$\alpha(M)$ to be
$$\alpha(M)=\liminf_{r\to\infty}\frac{i(p,r)}{r},$$
which is called the {\it injectivity radius growth} of 
$M$. Because in the definition of $\alpha(M)$ $r$ goes to 
infinity and the distance from $p$ to any other fixed point 
is a definite finite number, it can be seen directly 
(see also a proof in \cite{SW14}) that $\alpha(M)$ is 
not depending on $p$. One of their questions in \cite{SW14} 
is the following
\begin{ques}[\cite{SW14}]
For a complete noncompact manifold $M$, can one prove that 
every geodesic $\gamma:(-\infty, +\infty)\to M$ is a line 
as long as $\alpha(M)>1$?
\end{ques}

In other words, they asked that whether the injectivity 
radius of every point in $M$ is infinity when 
$\alpha(M)>1$? A positive answer of Question 1 directly 
follows from Proposition 2 below.

\begin{prop}
Let $M$ be a complete Riemannian manifold and $p\in M$. 
If for some $r>0$, $p$ is not a cut point of any point 
$x$ such that $\operatorname{d}(x,p)=r$, then the 
injectivity radius $i_p$ at $p$ $> r$.
\end{prop}
\begin{rem3}
The point in proving Proposition 2 is to show that the 
minimal geodesics for $p$ to points in the metric sphere 
$S_r(p)=\{x\in M: \operatorname{d}(p,x)=r\}$ covers the 
whole ball $B_r(p)=\{x\in M: \operatorname{d}(p,x)\le r\}$.
Though the conclusion of Proposition 2 may be already known 
by some experts, it seems that it is still not well-known 
and there is no proof can be found in the earlier 
literature. That is the reason why I decided to 
write down a proof.
\end{rem3}
\begin{rem4}
It can be proved that for $p\in M$ and $r>0$, if the 
minimizing geodesic from $p$ to each point $x$ such that 
$\operatorname{d}(x,p)=r$ is unique, then the injectivity 
radius of $p$ $\ge r$. However, the proof is more complicate 
than that of Proposition 2. So we will not go into that 
case here.
\end{rem4}

\begin{proof}[Proof of Proposition 2]
Let $T_p^1M$ denote the set of all unit vectors at $p$ in 
$M$.
Let us denote 
$$A(p,r)=\{X\in T_p^1M: \exp_p(tX)\text{ is minimal on 
$[0,r')$ for some $r'>r$} \}.$$
It suffices to show that $A(p,r)$ is open and close in 
$T_p^1M$. 
Firstly, it is well-known that the function 
$\sigma:T_p^1M\to \Bbb 
R^+$,
$$\sigma(X)=\sup\{t: \exp tX \text{ is minimal on 
$[0,t]$}\},$$
is continuous (see 2.1.5 Lemma in \cite{Klin82}). Hence by
definition $A(p,r)$ is open.

Now let us show that $A(p,r)$ is closed in $T_p^1M$. Assume 
a 
sequence of unit vectors $X_i\in A(p,r)$ converges to a 
unit 
vector $X\in T_p^1M$, 
then the geodesic $\exp(tX_i)$ converges to  $\exp(tX)$ 
point-wisely. Because all geodesic $\exp(tX_i)$ is minimal 
on 
$[0,r]$, the limit $\exp(tX)$ is also a minimal geodesic on 
$[0,r]$, 
and thus $\operatorname{d}(\exp(rX),p)=r$.
Moreover, by the assumption of Proposition 2, 
$\exp(rX)$ is not a cut point of $p$. Hence, there is 
$\epsilon>0$ 
such that $\exp(tX)$ is also minimal on $[0,r+\epsilon]$. 
Thus $X\in 
A(p,r)$ and $A(p,r)$ is closed.

Because $A(p,r)$ is open and closed, it coincides with 
$T_p^1M$. 
Therefore the injectivity radius at $p$ is $>r$.
\end{proof}

The following corollaries directly follows from Proposition 
2.
Recall that $p$ is called a pole if the injectivity radius 
of $p$ is infinity. In particular, $M$ is diffeomorphic to 
$\Bbb R^n$ by the exponential map $\exp_p:T_pM\to M$ at a 
pole.

\begin{cor5}
Let $M$ be a complete non-compact manifold. 
$M$ possesses a pole at $p$ if ( and only if )
there is a sequence $r_k\to \infty$ such that $p$ is 
not a cut point of any point in $S(p,r_k)$.
\end{cor5}

\begin{cor6}
$$\limsup_{r\to\infty}\frac{i(p,r)}{r}>1,$$ 
implies that every point in $M$ is a pole. Hence 
$\text{either }\limsup_{r\to\infty}\frac{i(p,r)}{r}\in 
[0,1],$ 
$\text{or } \limsup_{r\to\infty}\frac{i(p,r)}{r}=\infty.$
\end{cor6}

Because $\alpha(M)\le 
\limsup_{r\to\infty}\frac{i(p,r)}{r}$, Corollary 6 not only 
answers Question 1, but also strength the homeomorphism 
result of Theorem 1.2 in \cite{SW14} to diffeomorphism in 
the case of dimension 4.

\section*{Acknowledgement} 
The author thanks Zhongyang Sun for the discussion to 
bring the author's attention to the problem, and also 
thanks Jianming Wan for his interest that gave the author 
motivation to complete this note.
\bibliographystyle{plain}
\bibliography{growthradius}

\end{document}